\documentclass[10pt]{amsart}
\usepackage[]{amsmath, amsthm, amsfonts, amssymb}
\usepackage[all]{xy}

\newcommand {\N} {{\mathbb N}}
\newcommand {\C} {{\mathbb C}}
\newcommand {\R} {{\mathbb R}}
\newcommand {\RR} {{\mathcal R}}

\newcommand {\Z} {{\mathbb Z}}
\newcommand {\Q} {{\mathbb Q}}
\newcommand {\PP} {{\mathbb P}}

\newcommand {\E} {{\mathcal E}}
\newcommand {\I} {{\mathcal I}}
\newcommand {\F} {{\mathcal F}}
\newcommand {\G} {{\mathcal G}}
\newcommand {\Ci} {C^\infty}

\newcommand {\dt} {\bullet}
\newcommand {\HH} {{\mathbb H}}

\newcommand {\LL} {{\mathcal L}}

\newcommand {\BS} {{\mathbb S}}
\newcommand {\M} {{\mathcal M}}
\newcommand {\X} {{\mathcal X}}

\newtheorem{thm}[subsection]{Theorem}
\newtheorem{cor}[subsection]{Corollary}
\newtheorem{lemma}[subsection]{Lemma}
\newtheorem{prop}[subsection]{Proposition}

\newtheorem{remark}[subsection]{Remark}
\newtheorem{ex}[subsection]{Example}

\begin{document}

 \title{Partial Regularity and Amplitude}

\author{ Donu Arapura}
\address{Arapura: Department of Mathematics\\
Purdue University\\
West Lafayette, IN 47907}
\email{dvb@math.purdue.edu}
\thanks{Partially supported by the NSF}

\maketitle
\begin{abstract}
  The author continues the study of Frobenius amplitude, introduced
in an earlier paper, and compares this to various other positivity
notions. These are applied to deduce vanishing theorems.
A refinement of Castelnuovo-Mumford regularity is also given.
\end{abstract}

In \cite{arapura}, this author defined the {\em  Frobenius amplitude} $\phi(\E)$ of an
 algebraic vector bundle $\E$. This is 
 an integer which is defined by reduction modulo $p$ and  
which provides a measure of the positivity of $\E$.
This notion was introduced in order to formulate, and prove, the
following vanishing theorem:  on a smooth projective variety $X$ 
defined over a field of characteristic $0$,
 $H^q(X,\Omega_X^p\otimes \E)=0$ for $p+q> \dim X+\phi(\E)$. 
The main goal of this paper is to demonstrate the efficacy of this result 
(and its refinements) by extracting generalizations of a number of
known  theorems. The  point is to replace the Frobenius amplitude by a more
accessible expression involving the {\em amplitude} $\alpha(\E)$, which
is defined below. The basic result proved in section~\ref{section:vanish}, which
 subsumes the vanishing theorems of
Le Potier,  Sommese and others, is that the above groups vanish for $p+q>\dim X+
rank(\E)+\alpha(\E)$ when $\E$ satisfies some additional assumptions.

Over a field of characteristic $p>0$,  the Frobenius powers 
$\E^{(p^n)}$  of a vector bundle  $\E$ can be defined by raising
its transition matrices to $p^n$. The Frobenius amplitude $\phi(\E)$
is the threshold for which 
$H^i(\E^{(p^n)}\otimes \F) =0$ for any  coherent  $\F$, $n\gg 0$ and
$i> \phi(\E)$.
In general, the definition is indirect since it involves 
specialization into positive characteristic.  Replacing Frobenius powers
by symmetric powers in the above definition, leads to   the 
{ amplitude} $\alpha(\E)$ of $\E$. While this invariant is probably
new, it  is closely related  to a number of   preexisting positivity notions. 
For example, the amplitude of a vector bundle is at most $k$ 
if it is $k$-ample in Sommese's sense, or if  its dual is strongly $k+1$-convex in
 Andreotti and Grauert's sense when the field is $\C$. The technical
 heart of this paper is section \ref{section:estimates}, where  
the bounds given in  \cite[theorems 6.1, 6.7, 7.1]{arapura} are
refined by using  amplitude in the place of ampleness. In particular,
 the key inequality  $\phi(\E)\le rank(\E)+\alpha(\E)$ is obtained
 under appropriate assumptions.

The notion of   partial Castelnuovo-Mumford regularity
arose as a technical tool in the course of carrying out the above program.
However, it gradually became clear  that this concept is interesting
on its own. In rough terms, for partial regularity
 the conditions in Mumford's  definition
of $0$-regularity are required to hold for cohomology groups in
degrees greater than a number called the {\em level}
of the sheaf. 
The key point is that, unlike the notion of amplitude,
 the level involves a finite number of conditions, 
so it is easier to control in
families. The level is $0$ precisely for $0$-regular sheaves. 
Several well known results about regular
sheaves generalize quite nicely.
For example, the fact that the tensor product of regular vector bundles
on $\PP^n$ is regular, generalizes to the  subadditivity of the level
under  tensor products. 
The main ingredient in the proof of
this and related properties is the Beilinson-Orlov  (\cite{beil}, \cite{orlov})
resolution of the  diagonal.

In the final section, I consider the special class of varieties
 which admit endomorphisms of degree prime to the characteristic. 
 These include Abelian and toric varieties.
For such varieties  a vanishing theorem is obtained which is much
 stronger than anything available  for general varieties. 
The result generalizes known facts
about Abelian and toric varieties going back to Mumford and Danilov.
The proofs, however, have  an entirely different character; they  are
modeled on  Frobenius splitting arguments \cite{mehta}, with the endomorphism
playing the role of the Frobenius.

I would like to  thank  A. Dhillon for numerous discussions about varieties 
with endomorphisms, Y. Kawamata, D. Keeler and K. Matsuki
for various  helpful conversations and emails,  and  the referee for
catching a number of errors.

\section{Partial regularity}

Let $k$ be a field.
Recall that a coherent sheaf $\F$ on $\PP^n_k$ is {\em $0$-regular},
or simply just regular, if  $H^i(\F(-i)) = 0$ for all $i>0$ \cite{mumford1}.
We measure the deviation from regularity by defining the
 {\em level} of $\F$ to be
$$\lambda(\F) = 
 \max(\{q\mid \exists i\ge 0, H^{q+i}(\PP^n,\F(-1-i))\not=0\}\cup \{0\})$$
Equivalently, $\lambda(\F)$ is the smallest natural number for which $q>\lambda(\F)$
implies
$$H^q(\F(-1)) = H^{q+1}(\F(-2)) =\ldots =0$$
It follows that $\lambda(\F) =0$ if and only if $\F$ is regular.
An alternative  description of the level is that it measures the
complexity of Beilinson's spectral sequence (remark \ref{remark:beil}).
The point is that regularity is known to be  equivalent
to this spectral sequence collapsing to a single row concentrated
along the horizontal axis \cite{chips}.

\begin{lemma}\label{lemma:lambdaext}
If 
$$0\to\E_1\to \E_2\to \E_3\to 0$$ is an extension, then
$\lambda(\E_2) \le max(\lambda(\E_1),\lambda(\E_3))$.  
\end{lemma}

\begin{proof}
  This should be clear from the definition of $\lambda$ and the
exact sequence for cohomology.
\end{proof}

Let $X$ be a smooth projective variety defined over $k$.
Fix a very ample line bundle $\LL$ on $X$, and set $O_X(1) =
\LL^{\otimes N}$ 
for some $N>0$ to be specified later. The level of a sheaf on $X$
will be measured using $O_X(1)$.
Let $A= \oplus H^0(X,O_X(i))$
denote the projective coordinate ring.
The vector spaces $B_m$ are defined
inductively by $B_0=k$, $B_1= H^0(X,O_X(1))$ and
$$B_m = ker[B_{m-1}\otimes H^0(X,O_X(1))\to B_{m-2}\otimes
H^0(X,O_X(2))]$$
We can produce a complex of graded $A$-modules
\begin{equation}
  \label{eq:koszul}
 B_M\otimes A(-M+1)\to \ldots B_1\otimes A(-1)\to A\to k\to 0
\end{equation}

Set $\RR_0= O_X$, and
$$\RR_m = ker[B_m\otimes O_X\to B_{m-1}\otimes O_X(1)]$$
for $m>0$. 
From (\ref{eq:koszul}), we obtain a complex of sheaves
\begin{equation}
  \label{eq:koszul2}
0\to \RR_m\otimes_{O_X} O_X(-m+1)\to B_m\otimes_k O_X(-m+1)\to \ldots 
B_1\otimes_k O_X\to O_X(1)\to 0
\end{equation}

\begin{thm}[Orlov]\label{thm:orlov}
 Let $O_X(1) = \LL^{\otimes N}$ as above. Fix $M\ge m>0$.
 Then, there is a complex of sheaves on $X\times X$
\begin{equation}\label{eq:koszul3}
\RR_M\boxtimes O_X(-M)\to \ldots \RR_1\boxtimes O_X(-1)\to  \RR_0\boxtimes
O_X \to O_\Delta\to 0
\end{equation}
where $\Delta \subset X\times X$ is the diagonal.
The complexes (\ref{eq:koszul}),
(\ref{eq:koszul2}) and (\ref{eq:koszul3}) are exact for $N\gg 0$.
\end{thm}

\begin{proof}
Orlov \cite[prop A.1]{orlov} shows that 
the obstructions to exactness of the above complexes 
 lie in a finite number of coherent
cohomology groups, which vanish for large $N$ by Serre.
\end{proof}

\begin{remark}
  We are attributing this result to Orlov, since his paper is (as
far as we know) the first place where it has been stated in this form.
However there is some closely related  prior work that ought to
be mentioned. The ring $A$ has been shown to be  Koszul  for large $N$ by
 Backelin \cite{backelin}, and this implies the exactness of
(\ref{eq:koszul}) and (\ref{eq:koszul3}) for  all $M$ 
\cite[section 3]{kawamata}. The theorem can also be
deduced from the work of Ein and Lazarsfeld \cite{einL}.
\end{remark}

Standard semicontinuity  and generic flatness
arguments yield:

\begin{cor}\label{cor:orlov}
  Let  $\X\to T$ be a flat family over an irreducible quasicompact base,
and let $\LL$ be a relatively very ample line bundle on $\X$.
 For any $M\ge m>0$, there exists an $N$ and a nonempty open set $U\subset T$
such that for $O_\X(1) = \LL^{\otimes N}$ the complexes (\ref{eq:koszul}),
(\ref{eq:koszul2}) and (\ref{eq:koszul3}) are exact on all the
 fibers $\X_t$ over the closed points $t\in U$.
\end{cor}

We will say that $O_X(1)$ is {\em sufficiently ample} if
the complexes (\ref{eq:koszul}), (\ref{eq:koszul2}) and (\ref{eq:koszul3})
are for exact $M=\dim X$. Orlov's theorem guarantees that a high power 
of a fixed  ample bundle is sufficiently ample, and that this true
for all nearby fibers if $X$ varies in a family.

The fundamental example is the following.

\begin{ex}
  Let  $X =\PP^n$. The hyperplane bundle $O_X(1)$ is sufficiently ample. 
The coordinate ring $A\cong  k[x_0,\ldots x_n]$. The complex (\ref{eq:koszul})
is just the Koszul resolution of $A/(x_0,\ldots x_n)$,
  $\RR_m \cong \Omega^m(m)$, and  (\ref{eq:koszul3})
is  Beilinson's resolution  \cite{beil}:
 $$0\to \Omega^n(n)\boxtimes O(-n)\to \ldots \Omega^1(1)\boxtimes O(-1)\to
 O\boxtimes O$$
 of  $O_\Delta$. 
\end{ex}

This resolution is the basis for 
Beilinson's spectral sequences. One of which is:
\begin{equation}\label{eq:beil}
E_1^{ab} = H^b(\PP^n,\Omega^{-a}(-a)\otimes \E)\otimes O(a)\Rightarrow 
 \E
\end{equation}
Theorem \ref{thm:orlov} allows us to construct a similar spectral
sequence  for any $X$.

\begin{lemma}\label{lemma:tensor}
Suppose that $O_X(1)$ is sufficiently ample. Let
  $\E$ and $\F$ be  coherent sheaves with one of them  locally
  free. Then if for each pair of integers $a\ge 0$ and $b\ge 0$, one of
$H^b(\E\otimes \RR_a) = 0$
or 
$H^{i+a-b}(\F(-a)) = 0$,
 then 
$$H^i(\E\otimes \F) = 0$$
\end{lemma}

\begin{proof}
Tensoring Orlov's resolution (\ref{eq:koszul3})
 with $\E\boxtimes \F$ yields an exact sequence
$$\ldots (\E \otimes \RR_1)\boxtimes \F(-1)\to  (\E\otimes \RR_0)\boxtimes
\F \to \E\otimes \F\to 0$$
Breaking this into short exact sequences and writing out the long
exact sequences for cohomology (or by using a spectral sequence), we obtain
$H^i(\E\otimes \F) = 0$ provided that 
$$H^{i+a}(\E \otimes \RR_a\boxtimes \F(-a)) = 0$$
for all $a\ge 0$. The lemma now 
follows from K\"unneth's formula.
\end{proof}

We can control the cohomology groups in the preceding lemma by the following.

\begin{lemma}\label{lemma:lambdaRi}
Suppose that $O_X(1)$ is sufficiently ample.
For any coherent sheaf $\F$ and $a\ge 0$,  $H^j(X, \RR_i\otimes  \F(-a)) = 0$
for $j> \lambda(\F)+a$ and $0\le i \le \dim X$.
\end{lemma}

\begin{proof}
  We prove this by descending induction on $a$, starting from
$a=\dim X$ where it is trivial. The sequence (\ref{eq:koszul2})
yields
$$0\to \RR_{i+1}\otimes O_X(-1)\to B_{i+1}\otimes O_X(-1)\to 
\RR_i\to 0$$
This implies
$$ B_{i+1}\otimes H^{j}( \F(-a-1))\to
H^{j}(X, \RR_i\otimes \F(-a)) \to
H^{j+1}(X, \RR_{i+1}\otimes \F(-a-1))$$
and the lemma follows.
\end{proof}

\begin{cor}\label{cor:lambdaFa}
 $H^j(X,  \F(-a)) = 0$
for $j> \lambda(\F)+a$. 
\end{cor}

\begin{cor}\label{cor:lambdaRi}
 $H^j(X, \RR_i\otimes \F) = 0$
for $j> \lambda(\F)$. 
\end{cor}

\begin{remark}\label{remark:beil}
The last corollary   implies that there are no  nonzero rows
of  Beilinson's spectral sequence (\ref{eq:beil})
other than the  first $\lambda$ rows.
\end{remark}

  \begin{prop}\label{prop:lambdatensor}
    Suppose that  $O_X(1)$ is chosen sufficiently ample, and
 suppose that $\E$ and $\F$ are two coherent sheaves on $X$ with
 $\E$ locally free. Then $\lambda(\E\otimes \F) \le \lambda(\E)+\lambda(\F)$
  \end{prop}

  \begin{proof}
It suffices to prove that
\begin{equation}
  \label{eq:Hqi}
H^{q+i}(\E\otimes \F(-i-1)) = 0
\end{equation}
for $q>\lambda(\E)+\lambda(\F)$ and $i\ge 0$. Given $b\ge 0$, we
must have $b>\lambda(\E)$ or $q-b > \lambda(\F)$. In the first
case, $H^b(\E\otimes \RR_a) = 0$ by corollary \ref{cor:lambdaRi}.
In the second case, 
$$H^{q+i+a-b}(\F(-i-a-1)) = 0.$$
Therefore, lemma \ref{lemma:tensor} implies that (\ref{eq:Hqi}) holds. 
  \end{proof}

Recall that a coherent  sheaf $\F$ is called $m$-regular if $\F(m)$ is regular.
The  regularity $reg(\F)$ is the 
smallest $m$ such that $\F(m)$ is regular.
 We note that
$H^i(\F) = 0$ for $m\ge reg(\F)$ and $i>0$ \cite[p. 100]{mumford1}.

The following result is well known for $\PP^n$ \cite{lazarsfeld}.
  \begin{cor}
    Given $\E,\F$ and $O_X(1)$ as above.
If $\E$ is $p$-regular and
$\F$ is $q$-regular, then $\E\otimes \F$ is $p+q$-regular. 
  \end{cor}

  \begin{proof}
By assumption,
   $$\lambda(\E\otimes \F(p+q)) \le \lambda(\E(p))+\lambda(\F(q))=0.$$
  \end{proof}

\begin{lemma}\label{lemma:lambdatensorreg}
If $\E$ and $\F$ are as above. Then $\lambda(\E\otimes
\F)\le \lambda(\E(-reg(\F)))$
and $H^i(\E\otimes \F) = 0 $ for $i> \lambda(\E(-reg(\F)))$.
\end{lemma}

\begin{proof}
Let $r=reg(\F)$. Proposition~\ref{prop:lambdatensor} implies
$$\lambda(\E\otimes \F) =\lambda(\E(-r)\otimes \F(r)) \le
\lambda(\E(-r)) + 0$$
Corollary~\ref{cor:lambdaFa} shows that $H^i(\E\otimes \F) = 0 $ 
for $i> \lambda(\E(-r))\ge \lambda(\E\otimes \F)$
\end{proof}

\section{Amplitude}\label{section:amp}

Let $k$ be a field, and let $X$ be a
projective variety defined over $k$. Let  $Vect(r)=Vect(X,r)$  be the
set of isomorphism classes of vector bundles of rank $r$ over $X$,
and let $Vect(X) = \cup_r Vect(X,r)$.
Let $P^N:Vect(X)\to Vect(X), N\in \N$ be a sequence of operations. 
We define the {\em $P$-amplitude}
 $\alpha_{P}(\E)$ of a vector bundle $\E$ to be the 
 smallest integer $i_0$ such that for every locally free sheaf $\F$
there exists an $N_0$ such that $H^i(X,P^{N}(\E)\otimes {\F} ) = 0$
for $i> i_0$ and $N\ge N_0$. Of particular interest are 
the {\em amplitude}  $\alpha(\E)$ and the {\em $\otimes$-amplitude}
 $\alpha_\otimes(\E)$,
where the operations $P^\dt$ are the symmetric and tensor powers
respectively.
 If $char\, k = p>0$, we also have the Frobenius amplitude $\phi$ where
 $P^\dt$ are the Frobenius powers.
 In more explicit terms, these are the smallest integers for which
$$H^i(X,S^N(\E)\otimes \F) = 0, \mbox{ for }  N\gg 0,\, i> \alpha(\E)$$
$$H^i(X,\E^{\otimes N}\otimes \F) = 0, \mbox{ for } N\gg 0,\,  i> \alpha_\otimes(\E)
$$
and
$$H^i(X,\E^{(p^N)}\otimes \F) = 0, \mbox{ for } N\gg 0,\,  i> \phi(\E)$$
for every locally free sheaf $\F$.

The next example will be taken up again in the last section.
 Suppose that $f:X\to X$ is an endomorphism of schemes, then
we define the {\em $f$-amplitude} $\alpha_f(\E) = \alpha_P(\E)$ with
$P^n = (f^*)^n$. When $char k = p>0$, and $F:X\to X$ is the absolute
Frobenius, then $\alpha_F$ coincides with the Frobenius amplitude
because $(F^{*})^n\E = \E^{(p^n)}$.

We list a few basic properties of these invariants.

 \begin{lemma}\label{lemma:alphaOE}
$\alpha(\E) = \alpha(O_{\PP(\E)}(1)) =\alpha_\otimes(O_{\PP(\E)}(1))$.
\end{lemma}

\begin{remark}
We are using the 
convention $\PP(\E) = {\bf Proj}(S^*(\E))$.  
\end{remark}

\begin{proof}
We have $\alpha(O_{\PP(\E)}(1)) =\alpha_\otimes(O_{\PP(\E)}(1))$
since the symmetric and tensor power powers of $O_{\PP(E)}(1)$
coincide. The isomorphism
\begin{equation}
  \label{eq:cohomologyofspecial}
H^i(X, S^N(\E)\otimes \F)  \cong H^i(\PP(\E), O_\PP(N)\otimes
\pi^*\F)= H^i(\PP(\E), O_\PP(1)^{\otimes N} \otimes
\pi^*\F)  
\end{equation}
where $\pi:\PP(\E)\to X$ is the projection, implies $\alpha(\E) \le
\alpha(O_{\PP(\E)}(1))$. For  the opposite inequality, it suffices to
prove that $H^i(\PP(E), O_\PP(N)\otimes \mathcal{G})=0$ for
all coherent $\mathcal{G}$ on
$\PP(E)$, $i> \alpha(\E)$ and $N\gg 0$. We do this by descending $i$.
Call a sheaf {\em special} if  it isomorphic to a sheaf of the
form $O_{\PP}(j)\otimes \pi^*\F$. For  special sheaves, the desired
vanishing statement follows from (\ref{eq:cohomologyofspecial}).
In general, observe that for any coherent sheaf
$\mathcal{G}$ on
$\PP(E)$, we have a surjection  $\pi^*\pi_*\mathcal{G}(M)\to
\mathcal{G}(M)$ for $M\gg 0$.
This is  just a relative form of Serre's global generation
theorem.  Thus we have an exact sequence
$$0\to \mathcal{G}'\to \mathcal{G}''\to \mathcal{G}\to 0$$
where the middle sheaf is special. The vanishing of  $H^i(
O_\PP(N)\otimes \mathcal{G})$ for $i> \alpha(\E)$ now follows
 by descending induction.
\end{proof}

 \begin{lemma}\label{lemma:alphaO1}
   Fix an ample line bundle $O_X(1)$. Then the following are
equivalent:
\begin{enumerate}
\item[(a)] $\alpha_{P}(\E)\le A$.
\item[(b)]  For any $b$ there exists $N_0$ such that 
$$H^i(P^N(\E)(b)) = 0$$
for all $i> A$, $N\ge N_0$.
\item[(c)] For any coherent sheaf $\F$ there exists $N_0$ such that 
 $$H^i(X,P^{N}(\E)\otimes {\F} ) = 0$$
for all $i> A$, $N\ge N_0$.
\end{enumerate}
\end{lemma}

\begin{proof}
  To prove that (b) implies (c), we use the fact that any coherent sheaf
$\F$ can be resolved as
$$0\to \F_1\to \bigoplus_i O_X(b_i)\to \F\to 0$$
Statement (c) can then be proven by descending induction on $A$. 
The other implications are automatic. 
\end{proof}
 
\begin{lemma} If $\E$ is a vector bundle, then $\alpha(\E) =0$ 
if and only if $\E$ is ample. If $char\, k=0$, then
$\alpha_{\otimes}(\E)=0$ if and only if $\E$ is ample.
\end{lemma}

\begin{proof} This follows from \cite[prop 3.3, section 5]{hartshorne}
    and the previous lemma.
 \end{proof}  

 \begin{lemma}\label{lemma:Tnample}
If $\E$ is a vector bundle and $n$ a positive integer, then
   $$\alpha_\otimes(\E^{\otimes n})= \alpha_\otimes(\E).$$
 \end{lemma}

 \begin{proof}
  Let $\F$ be a vector bundle. Then we can choose $N_0>0$
such that 
$$H^i(X, \E^{\otimes nN+a}\otimes \F)=
H^i(X, (\E^{\otimes n})^{\otimes N}
\otimes (\E^{\otimes a}\otimes \F)) = 0$$
for $i>\alpha_\otimes(\E^{\otimes n})$, $N>N_0$, $0\le a< n$.
This shows that $\alpha_\otimes(\E^{\otimes n})\ge
\alpha_\otimes(\E)$. For the opposite inequality, observe
that 
$$H^i( (\E^{\otimes n})^{\otimes N}\otimes \F) =0$$
for $i> \alpha(\E)$ and $N\gg 0 $.
 \end{proof}

 \begin{cor}\label{cor:Snample}
 If $\E$ is a vector bundle and $n$ a positive integer, then
   $\alpha(S^n(\E))\ge \alpha(\E)$.  
 \end{cor}

 \begin{proof}
   This follows by applying the lemma to $O_{\PP(\E)}(1)$ and
   appealing to lemma~\ref{lemma:alphaOE}.
 \end{proof}

The following is a generalization of \cite[2.4(4)]{arapura}:
 
\begin{lemma}\label{lemma:alphaPfE}
 Let $f:X\to Y$ be a morphism of projective varieties.
  Suppose that $P^\dt:Vect(X)\to Vect(X)$ and $\tilde P^\dt:Vect(Y)\to 
  Vect(Y)$ be a collection of operations which are compatible in
  the sense that $P^\dt(f^*\E) = f^*\tilde P^\dt (\E)$. Then
  $$\alpha_{P}(f^*\E) \le  \alpha_{\tilde P}(\E) + d$$
  where $d$ is the maximum of the dimensions of the fibers.
\end{lemma}  

\begin{proof}
 Let $\F$ be a locally free $O_X$-module. We use
the  Leray spectral sequence 
$$E_2^{ij}= H^i(Y, \tilde P^N(\E)\otimes R^jf_*\F)
\Rightarrow H^{i+j}(X,P^N(f^*\E)\otimes \F) $$
By assumption $E_2^{ij}=0$ when $j>d$. Also
$E_2^{ij}=0$ for  $i>\alpha_{\tilde P}(\E)$, all $j$
 and $N \gg0$, by lemma \ref{lemma:alphaO1}.
 Therefore 
 $$ H^{i}(X,P^N(f^*\E)\otimes \F)  =0$$
 for $i>\alpha_{\tilde P}(\E)+d$ and $N \gg 0$.
 \end{proof} 

 \begin{cor}\label{cor:restalpha}
   If $X\to Y$ is a closed immersion, then 
$\alpha(\E|_X)\le \alpha(\E)$.
 \end{cor}

 \begin{cor}\label{cor:quot}
   If $\E\to \F$ is a surjection of vector bundles, then
$\alpha(\F)\le \alpha(\E)$
 \end{cor}

 \begin{proof}
We have an inclusion $\PP(\F)\subset \PP(\E)$ such that
$$O_{\PP(\E)}(1)|_{\PP(\F)} = O_{\PP(\F)}(1).$$   
Now apply the previous corollary.
 \end{proof}

There is a log version of amplitude, analogous to the notion
of  Frobenius amplitude  $\phi(\E,D)$
of $\E$ relative to $D$ given in \cite{arapura}.
 Let $X$ be a smooth projective 
variety with a reduced normal crossing divisor $D$. Given a vector 
bundle $\E$ on $X$, define the {\em amplitude of $\E$ relative to $D$} by
$$\alpha(\E, D) = 
\min \{ \alpha(S^n(\E)(-D'))\mid n\in \N, 0 \le D'\le (n-1) D\}$$
It is perhaps more instructive to view this as the minimum of the
amplitudes of the ``vector bundles''
$\E\langle -\Delta\rangle$  \cite[6.2A]{lazarsfeld}  as $\Delta$ ranges over strictly fractional
effective $\Q$-divisors with support in $D$.

\begin{lemma}\label{lemma:alphaEDonPE}
  Let $\pi:\PP(\E)\to X$ be the canonical projection, then
$\alpha(\E,D)= \alpha(O_{\PP(\E)}(1),\pi^*D)$.
\end{lemma}

\begin{proof}
 This is immediate from the definition. 
\end{proof}

\begin{lemma}\label{lemma:alphaEDbirat}
 Let $Y\to X$ be a morphism of projective varieties with $Y$ smooth.
Suppose that $D= \sum D_i$ is a divisor with normal
crossings on $Y$, such that there exist $a_i\ge 0$
for which $L=O_Y(-\sum a_iD_i)$ is
relatively ample. Then for any locally free sheaf $\E$
on $X$, $\alpha(f^*\E, D)\le \alpha(\E)$. 
\end{lemma}

\begin{proof}
 After replacing $f$ by $\PP(f^*\E)\to \PP(\E)$ and invoking the previous
lemma, we can assume that $\E$ is a line bundle. Choose $n_0>a_i$,
and set $\M = f^*\E^{\otimes n_0}\otimes L$. Let $\F$ be
a coherent sheaf on $Y$.
Since $L$ is relatively ample, the higher direct images
of $\F\otimes L^{\otimes n}$ vanish for $n\gg0$.
The spectral sequence
$$H^a(\E^{\otimes nn_0}\otimes R^bf_*(\F\otimes L^{\otimes n}))
\Rightarrow H^{i}(\M^{\otimes n}\otimes \F)$$
yields the vanishing of the abutment for $i>\alpha(\E)$ and
$n\gg 0$.
\end{proof}

We define the generic amplitude of a vector bundle
$\E$ on a projective variety $Y$ by
$$\alpha_{gen}(\E) =\inf \alpha(r^*\E, Ex(r)), $$
where $r:Y'\to Y$ varies over birational maps from smooth varieties resolutions with normal crossing
exceptional divisor $Ex(r)$. In contrast to $\alpha(\E)$,
this is a birational invariant.

\begin{lemma}\label{lemma:alphagen}
   Suppose that $char\, k=0$, and let $f:X\to Y$ be a morphism of 
projective varieties. Then
  $$\alpha_{gen}(f^*\E) \le  \alpha(\E) + d$$
 where $d$ is the dimension of the generic fiber.
\end{lemma}

\begin{proof}
 This can be deduced  from  lemma~\ref{lemma:alphaEDbirat} in exactly
the same way that corollary 2.8 was deduced from lemma 2.6 in
 \cite{arapura}.
\end{proof}

Recall that a vector bundle $\E$ is {\em $d$-ample} \cite{sommese} if 
some $O_{\PP(\E)}(m)$ is   base point free and dimensions of the
fibers of the induced map
    $$\phi_m:\PP(\E)\to \PP(H^0(O_{\PP(\E)}(m)))$$
 are all less than or equal to $d$. An immediate consequence of
the preceding lemma is:

 \begin{prop}[Sommese]\label{prop:sommese}
Suppose that some positive symmetric power $S^m(\E)$ is globally generated.
Then $\alpha(\E) \le  d$ if and only if $\E$ is $d$-ample.
\end{prop}

\begin{proof}
This is  proven in  \cite[prop 1.7]{sommese} under the blanket assumption that
$k =\C$. For any field, when $\E$ is $d$-ample,
 the inequality $\alpha(\E) \le  d$ 
is a consequence of  lemma \ref{lemma:alphaPfE} applied to 
$P=\PP(\E)$.
Conversely, suppose that $S^m(\E)$ and therefore that $O_{\PP(\E)}(m)$ is
globally generated. Then for any coherent sheaf $\F$ supported on a closed
fiber of $\phi_m$,
$$H^i(P,\F) \cong H^i(P,\F\otimes O_{\PP(\E)}(mN)) = 0$$
for $i>\alpha(\E)$ and $N\gg 0$. This implies that the dimension of the
fiber, which coincides with its coherent cohomological
dimension, is at most $\alpha(\E)$.

\end{proof}

Since the global generation assumption for  $S^m(\E)$ is built into the
definition of $d$-ampleness, the  above equivalence fails without it.
We want to refine this to get estimates of amplitude under a weaker assumption.
Consider the condition $G(\E,m,U)$ that for some $m$,
the map $H^0(X,S^m(\E))\otimes O_U\to S^m(\E)|_U$ is surjective
for some open set $U\subseteq X$. Under this assumption, we
get canonical morphisms
$$\phi_{m,U}:\PP(\E|_U)\to \PP(H^0(X,S^m(\E)))$$
and   
$$\psi_{m,U}:U\to Grass(H^0(S^m(\E)),rank(S^m(\E)) ),$$
where $Grass(V,r)$ is the Grassmanian of $r$-dimensional quotients
of $V$. These maps are related. If $Q$ denotes the universal quotient
bundle on the Grassmanian, $\psi_{m,U}$  induces a map
$\PP(\E|_U)\to \PP(Q)$, and $\phi_{m,U}$ is composite of
this and the projection $\PP(Q)\to \PP(H^0(X,S^m(\E)))$.

\begin{prop}\label{prop:gensommese}
  Suppose that $char\, k=0$ and that $\E$ is a vector bundle
on a smooth projective variety $X$ such that $G(\E,m,U)$ holds.
\begin{enumerate}
\item If the restriction of the universal bundle $Q$ to 
$\overline{\psi_{m,U}(U)}$ is $d$-ample, then
$$\alpha_{gen}(\E) \le d +  \dim X - \dim \psi_{m,U}(U)$$ 

\item If $\phi_{m,U}$ extends to a map $\overline{\phi}_{m,U}:Y\to\PP(H^0(X,S^m(\E))) $
of some nonsingular model of $\PP(\E)$ such that
the dimensions of  all the fibers of $\overline{\phi}_{m,U}$ are at most $d$,
then
$$\alpha(\E)\le max(d, (rank(\E)-1)+dim(X-U)).$$

\end{enumerate}
\end{prop}

We postpone the proof until after the next lemma.

 \begin{lemma}\label{lemma:basept}
 Suppose $char\, k =0$.
Let $\LL$ be a line bundle over a projective variety  $X$ with
at worst rational Gorenstein  singularities. Suppose that
$f:Y\to X$ is a desingularization with a globally generated sub-line
bundle $\M\subset f^*\LL$, and let $B\subset Y$ be the support of $\LL/f_*(\M)$.
Then 
$$\alpha(\LL) \le max(\alpha(\M), \dim B)$$
 \end{lemma}

 \begin{proof}

Let $\E$ be a locally free sheaf on $X$. Set $\E' = \E\otimes
\omega_X^{-1}$.
 By the Grauert-Riemenschneider vanishing theorem \cite{grauert}
 \begin{equation}
    \label{eq:omxM}
      \R f_* (\omega_Y \otimes \M^{\otimes n}) = 
f_*(\omega_Y\otimes \M^{\otimes n})
  \end{equation}
This  along with the projection formula yields a diagram with
an exact rows
$$
\xymatrix{
   H^i(\omega_Y\otimes \M^{\otimes n} \otimes f^*\E')\ar@{=}[d] &  &  \\
   H^i(\R f_*(\omega_Y\otimes \M^{\otimes n} \otimes
   f^*\E'))\ar@{=}[d] &  &  \\
   H^i(f_*(\omega_Y\otimes \M^{\otimes n}) \otimes \E')\ar[r]^r & H^i(\omega_X\otimes \LL^{\otimes n} \otimes \E')\ar[r] & H^i(Q\otimes \E')
}
$$
where $Q=\omega_X\otimes \LL^{\otimes n}/f_*(\omega_Y\otimes \M^{\otimes
  n})$.
Since $Q$ has  support in $B$,
 $r$ is surjective for $i>\dim B$, which implies 
$$H^i(\LL^{\otimes n} \otimes \E) = 0$$
for $n>>0$ and $i> max(\alpha(\M),\dim B)$.
 \end{proof}

 \begin{proof}[Proof of proposition \ref{prop:gensommese}]
The first inequality is a consequence of lemma~\ref{lemma:alphagen} and 
proposition~\ref{prop:sommese}. The second follows by applying
lemmas \ref{lemma:alphaOE}, \ref{lemma:Tnample}, and \ref{lemma:basept}
to $O_{\PP(\E)}(1)$.
 \end{proof}
 
\section{Subadditivity of amplitude}
We continue the assumptions from the previous section.

\begin{thm}\label{thm:ext}
If 
$$0\to\E_1\to \E_2\to \E_3\to 0$$
is an exact sequence of vector bundles, 
then $\alpha(\E_2)\le \alpha(\E_1)+\alpha(\E_3)$.   
 \end{thm}

 \begin{proof}
Fix  a sufficiently ample line bundle $O_X(1)$.
Let $\F$ be a vector bundle. 
Choose $M \gg 0$ so that $\F(M)$ is regular, i.e. 
\begin{equation}
  \label{eq:ext1}
\lambda( \F(M)) =0.  
\end{equation}
From the definition of $\alpha$ and $\lambda$, observe that given a bundle
$\E$ and an integer $m$, there exists an integer $n_0$ such that for $n\ge n_0$,
$$\lambda(S^n(\E)(m))\le \alpha(\E)$$ 
Therefore, we can find $p_0$ so that 
\begin{equation}
  \label{eq:ext2}
\lambda(S^p(\E_1)(-M))\le \alpha(\E_1)  
\end{equation}
\begin{equation}
  \label{eq:ext3}
\lambda(S^p(\E_3)) \le \alpha(\E_3)  
\end{equation}
for $p\ge p_0$.
Choose $N\gg 0$ so that 
\begin{equation}
  \label{eq:ext4}
 \lambda(S^p(\E_i)\otimes \F(N) )=0
\end{equation}
for $p<p_0$ and $i=1,3$.
Finally choose $r_0> p_0$ so that for $r\ge r_0$, 
\begin{equation}
  \label{eq:ext5}
 \lambda(S^{r-p}(\E_i)(-N)) \le \alpha(\E_i)  
\end{equation}
when $p < p_0$ and $i=1,3$. We claim that for $r\ge r_0$
$${\mathcal G}=S^p(\E_1)\otimes S^{r-p}(\E_3)\otimes \F$$
has level at most $  \alpha(\E_1)+\alpha(\E_3)$.
This follows by 
applying  proposition \ref{prop:lambdatensor}
to a decomposition $\G = \G_1\otimes \G_2\otimes \ldots$, to obtain
$$\lambda(\G) \le \sum \lambda(\G_i)\le \alpha(\E_1)+\alpha(\E_3) ,$$
in the three cases:
\begin{enumerate}
\item[I.] If $p< p_0$, then take 
$$\G_1=  S^{r-p}(\E_3)(-N),\> \G_2= S^{p}(\E_1)\otimes \F(N)$$
Use (\ref{eq:ext4}) and (\ref{eq:ext5}).
\item[II.] If $r-p< p_0$ then take
$$\G_1= S^p(\E_1)(-N), \> \G_2= S^{r-p}(\E_3)\otimes \F(N)$$
Use (\ref{eq:ext4}) and (\ref{eq:ext5}).
\item[III.] If $p, r-p \ge p_0$, then take
$$\G_1= S^p(\E_1)(-M), \> \G_2 =  S^{r-p}(\E_3), \> \G_3= \F(M)$$
Use (\ref{eq:ext1}),~(\ref{eq:ext2}) and (\ref{eq:ext3}).
\end{enumerate}

By \cite[III, exer. 5.16]{hart-book}, there is a filtration
$F^\dt \subseteq S^r(\E_2)$ such that
$$Gr^p(S^r(\E_2)) \cong S^p(\E_1)\otimes S^{r-p}(\E_3)$$
Repeated application of lemma \ref{lemma:lambdaext} and the
above claim will imply that $S^r(\E_2)\otimes \F$ will have
level at most $\alpha(\E_1)+\alpha(\E_3)$ for $r>r_0$. Lemma
\ref{lemma:lambdaRi} will imply that $H^i(S^r(\E_2)\otimes \F)=0$
for $i >\alpha(\E_1)+\alpha(\E_3)$ as required.
 \end{proof}

The following theorem is proved by the same method as above.

 \begin{thm}
Given vector bundles $\E_1$ and $\E_3$, we have
 $$\alpha_\otimes(\E_1\oplus \E_3)\le 
\alpha_\otimes(\E_1)+\alpha_\otimes(\E_3)$$
  $$\alpha_\otimes(\E_1\otimes \E_3)\le \alpha_\otimes(\E_3)
+\alpha_\otimes(\E_3)$$
\end{thm}

\begin{proof}
Set $\E_2 = \E_1\oplus \E_3$.
  The proof that 
$\alpha_\otimes (\E_2)\le \alpha_\otimes(\E_1)+\alpha_\otimes(\E_3)$
will follow the same outline as the proof of theorem~\ref{thm:ext}.
Let $\F$ be a vector bundle.
Choose $M \gg 0$ so that $\F(M)$ is regular, i.e. $\lambda( \F(M)) =0$.
We can find $p_0$ so that 
$$\lambda(\E_1^{\otimes p}(-M))\le \alpha_\otimes(\E_1)$$
$$\lambda(\E_3^{\otimes p}) \le \alpha_\otimes(\E_3)$$
for $p\ge p_0$.
Choose $N\gg 0$ so that the sheaves
 $\E_i^{\otimes p}\otimes \F(N)$ are regular
for $p<p_0$ and $i=1,3$.
Finally choose $r_0> p_0$ so that for $r\ge r_0$, 
 $\lambda(\E_i^{\otimes (r-p)}(-N)) \le \alpha_\otimes(\E_i)$
for $p < p_0$. We can then argue that for $r\ge r_0$
$${\mathcal G}=\E_1^{\otimes p}\otimes \E_3^{\otimes (r- p)}\otimes \F$$
has level at most 
$  \alpha_\otimes(\E_1)+\alpha_\otimes(\E_3)$, by splitting it
up into cases as above.
Lemma \ref{lemma:lambdaRi} will then imply that for 
 $ i> \alpha_\otimes(\E_1)+\alpha_\otimes(\E_3)$, 
$H^i(\G) =0$.  Therefore $H^i(\E_2^{\otimes r}\otimes \F)=0$,
since  $\E_2$ can be decomposed into a direct sum of sheaves
of the form $\E_1^{\otimes p}\otimes \E_3^{\otimes r- p}$.

By the first inequality 
$$H^i((\E_1\otimes \E_3)^{\otimes N}\otimes \F)
\subset H^i(\E_2^{\otimes 2N}\otimes \F) =0$$
for $ i> \alpha_\otimes(\E_1)+\alpha_\otimes(\E_3)$ and
$N\gg 0$. This implies the second inequality
$\alpha_\otimes(\E_1\otimes \E_3)\le \alpha_\otimes(\E_3)
+\alpha_\otimes(\E_3)$.
\end{proof}

\begin{prop}\label{prop:alphaVsalphaOtimes}
If $char\, k =0$, then for any vector bundle of rank $r$, we have
 $$\alpha(\E)\le \alpha_{\otimes}(\E)\le \alpha(\E^{\oplus r!}) $$
Equality holds if $\E$ is globally generated.  
\end{prop}

\begin{proof}
  Since $S^n(\E)$ is a direct summand of $\E^{\otimes n}$,
the inequality  $\alpha(\E)\le \alpha_{\otimes}(\E)$ holds.
By standard representation theoretic machinery \cite{FH}, the
 tensor power  $\E^{\otimes n}$ can be decomposed into a
direct sum of   Schur powers $\BS^\lambda(\E)$ for various
partitions $\lambda$ of $n$. Each $\BS^\lambda(\E)$
is a direct summand of some
$$S^{q_1}(\E)\otimes\ldots S^{q_{r!}}(\E), $$
and therefore of
$$ S^{q_1+\dots q_{r_!}}(\E^{\oplus r!}),$$
by \cite[prop 5.1]{hartshorne}, with $\sum q_i =n$. Thus 
$$H^i(S^{q_1}(\E)\otimes\ldots S^{q_{r!}}(\E)\otimes \F) = 
0$$
for $i> \alpha(\E^{\oplus r!})$ and $n\gg 0$.
Therefore
$$H^i(\E^{\otimes n}\otimes \F) = 0$$
for $i> \alpha(\E^{\oplus r!})$ and $n\gg 0$, as required.

The last part follows from \cite[cor 1.10]{sommese} and 
proposition \ref{prop:sommese}
\end{proof}

\begin{cor}
  $\alpha_\otimes(\E)\le r!\alpha(\E)$.
\end{cor}

\section{$q$-convexity}

In this section, we work over $\C$. By the GAGA theorem \cite{serre}, 
over a compact base
variety, we can replace algebraic vector bundles by the corresponding
analytic bundles. We will  denote a geometric vector bundle 
$\pi:E\to X$ with the standard font, and write $\E$ for its sheaf 
of holomorphic sections.

A complex manifold $Z$ is called {\em strongly
$q$-convex} (in the sense of Andreotti-Grauert) 
if  there is a $\Ci$ exhaustion $\psi:E\to \R$
such that the Levi form, given locally by
$$ L(\psi) = 
\left(\frac{\partial^2\psi}{\partial z_i\partial \overline{z_j}}\right), $$
has at most $q-1$ nonpositive eigenvalues at all points $z\in Z $
outside a compact  set. 
If $\E$ is a  holomorphic vector bundle,
we say that it is strongly $q$-convex if its total space  $E$ is.
A useful criterion is the following:

\begin{lemma}
  Let $\E$ be a holomorphic vector bundle with a Hermitean metric 
over a compact complex manifold $X$. Let $\Theta$ be the curvature  of the
compatible connection viewed as a $\Ci$ section of 
$E^*\otimes \overline{E}^*\otimes T_X^*\otimes \overline{T_X}^*$.
If  $\Theta(\xi,\bar \xi,-,-)$ has fewer than $q$ nonpositive  eigenvalues
for all $\xi\not= 0$ and all $x\in X$. Then $\E^*$ is strongly $q$-convex.
\end{lemma}

\begin{proof} The argument is indicated in \cite[p. 257]{andreotti}.
Nevertheless, we outline it here since our  notation is slightly different.
Since the signs of the curvatures of $\E$ and $\E^*$ are opposite,
 it suffices to  prove the dual statement that if  
$\Theta(\xi,\bar \xi,-,-)$ has
fewer than $q$ nonnegative eigenvalues then $\E$  is $q$-convex with respect
to $\psi$, where $\psi:E\to \R$ is the square of the norm.
  
We can choose local coordinates $z_i$  about $x$ and a local frame $e_\alpha$
of $\E$ such that the metric is given by matrix $h_{\alpha\beta}(z_1,\ldots)$
 satisfying
$$\frac{\partial h_{\alpha\beta}}{\partial z_i}|_x
 = \frac{\partial h_{\alpha\beta}}{\partial \bar z_i}|_x = 0.$$
In these coordinates, the curvature tensor at $(x,\xi)\in E$ is given by
$$ \Theta_{\alpha\beta ij} = 
-\frac{\partial^2 h_{\alpha\beta}}{\partial z_i\partial \bar z_j}$$
\cite[p. 119]{ss}, and
$$\psi = \sum h_{\alpha\beta}e_\alpha \bar e_\beta$$
The Levi form $L(\psi)$ at $(x,\xi)$ is represented by the matrix
$$ \begin{pmatrix} 
-\sum \Theta_{\alpha\beta ij}\xi_\alpha\bar \xi_\beta & 0\\ 
0& h_{\alpha\beta} \end{pmatrix} $$ 
and the lemma follows.
\end{proof}

\begin{thm}[Grauert]
  If $\E^*$ is a strongly $q$-convex vector bundle over a smooth complex 
projective variety, then $\alpha(\E)\le q-1$.
\end{thm}

The result follows from \cite[theorem 4.3]{grauertq}, which was
stated without proof, so we supply one here.
First, we observe:

\begin{prop}
  Suppose that $\pi:E\to X$ is a geometric holomorphic vector bundle over a 
compact complex manifold $X$, with associated locally free
sheaf $\E$. Then for any line bundle $\LL$ on $X$,
there exists a filtration, indexed by $\N$, on $H^i(E,\pi^*\LL)$ such that
$$Gr^\ell H^i(E,\pi^*\LL) \cong  H^i(X,S^\ell(\E^*)\otimes \LL)$$
\end{prop}

\begin{proof}
The argument  is a straightforward modification of the proof
of \cite[prop 26]{andreotti}.
 Choose an open  Leray 
 covering $\{U_i\}$ of $X$ which trivializes $\E$ and $\LL$.
 Then a section $\sigma$ of
  $\pi^*\LL$ over $\pi^{-1}U_i$ can be written as an infinite series
 $$\sum   \xi_1^{n_1}\ldots \xi_r^{n_r}\otimes \lambda_{n_1\ldots n_r} $$
 where $ \lambda_{n_1\ldots n_r}$ is a  section of $\LL(U_i)$ and
  $\xi_j$ are coordinates along the fibers of $E$. We filter this
 so that $\sigma\in F^\ell$ if  $ \lambda_{n_1\ldots n_r}=0$ for
 $\sum n_i< \ell$. We have a splitting
\begin{equation}\label{eq:splitpiLL}
\pi_*\pi^*\LL(U_i)=\pi^*\LL(\pi^{-1} U_i) = P_\ell(\pi^{-1} U_i)\oplus F^\ell(\pi^{-1} U_i) 
\end{equation}
where $P_\ell(\pi^{-1} U_i)$ is the space of polynomials (i.e. finite sums)
of degree at most $\ell$. Therefore $F^\ell/F^{\ell+1}$ can be identified
with the space of homogeneous polynomials of degree $\ell$.

In more invariant language, $F^\dt$ is the $\I$-adic filtration, where
$\I$ is the ideal of the zero section of $E$. In particular, the filtration
globalizes. We can give $H^i(\pi^*\LL)$ the induced filtration.  
Since the higher direct images $R^i\pi_*$, $i>0$, vanish for coherent sheaves, we can
identify  $H^i(\pi^*\LL)\cong H^i(\pi_*\pi^*\LL)$ with their
filtrations.
Then we have isomorphisms 
\begin{equation}\label{eq:GrLL}
Gr^\ell(\pi_*\pi^*\LL)\cong \pi_*Gr^\ell(\pi^*\LL) \cong \pi_*[(\I^\ell/\I^{\ell+1})\otimes \pi^*\LL] 
\cong S^\ell(\E^*)\otimes \LL
\end{equation}
The local splittings (\ref{eq:splitpiLL}) patch to give a global
splitting of the image of the sequence
$$ 0\to \I^\ell\pi^*\LL\to \pi^*\LL\to  \pi^*\LL/\I^\ell\pi^*\LL 
\to 0$$
under $\pi_*$.
Therefore, we have a diagram with exact rows
$$\xymatrix{
  & H^i(\I^\ell\pi^*\LL)\ar[r]\ar[d]^{\cong} &
 H^i(\pi^*\LL)\ar[r]\ar[d]^{\cong} & H^i( \pi^*\LL/\I^\ell\pi^*\LL)
 \ar[d]^{\cong} & \\ 
 0\ar[r]&  H^i(\pi_*\I^\ell\pi^*\LL)\ar[r] & H^i(\pi_*\pi^*\LL)\ar[r] &
 H^i( \pi_*\pi^*\LL/\I^\ell\pi^*\LL) \ar[r] &0
}
$$
The proposition is a now consequence of (\ref{eq:GrLL}).
\end{proof}

\begin{proof}[Proof of theorem]
By the proposition,
$$\dim H^i(E,\pi^*O(j)) \ge \sum_{\ell=0}^L  H^i(X,S^\ell(\E)\otimes
O(j)),\quad \forall L$$
The left side is finite dimensional for $i\ge q$ and arbitrary $j$ by 
\cite[theorem 14]{andreotti}. Therefore $H^i(S^\ell(\E)\otimes O(j))=0$
for large $\ell$ and $i\ge q$. This implies the theorem by 
lemma \ref{lemma:alphaO1}.
\end{proof}

\section{ Estimates on Frobenius amplitude}\label{section:estimates}

 In positive characteristic, the   Frobenius amplitude, or more briefly
$F$-amplitude, of a bundle $\E$ on $X$ was defined at the beginning of 
section \ref{section:amp}.
Over a field $k$ of  characteristic $0$, the definition is via specialization; 
see \cite{arapura} for a detailed discussion.
Given a pair $(X, \E)$, a thickening of it is a flat family
 $(\tilde X,  \tilde \E)$ over the spectrum of a subring $A\subset k$
of finite type, such that the original pair is given by 
base change. Then the $F$-amplitude is defined so that
$\phi(\E)\le C$ if and only if there exists a thickening such that
$\phi(\tilde \E |_{X_q})\le C$ for all closed points $q\in Spec\, A$.
In the sequel, we will of often write $\E_q$ for $\tilde \E|_{X_q}$,
and $p(q)$ for the characteristic of $A/q$.

We recall the notion of {\em  $F$-semipositivity} of vector bundles
from \cite[3.4]{arapura}.
If $char\, k=p>0$, $\E$ is $F$-semipositive if
the regularities of $\E^{(p^n)}$ are bounded away from $+\infty$.  
If $char\, k=0$, $\E$ is $F$-semipositive if there exists a thickening 
$(\tilde X,\tilde \E)$ such that  the restrictions $\tilde \E|_{X_q}$
are $F$-semipositive.
The definition is a technical one. The 
 key properties to keep in mind are 
that the notion is stable under pullbacks,
and tensor product with a $F$-semipositive bundle does not increase
$\phi$ [loc. cit.].

Moreover, a line bundle $\LL$ is $F$-semipositive if and only
if it is {\em arithmetically nef}, which means that either $char\, k>0$ and
$\LL$ is nef in the usual sense, or there is a thickening
$(\tilde X,\tilde \LL)$ such that all specializations  $\LL_q$ are nef
\cite[3.12]{arapura}, \cite{keeler2}. The proof of equivalence is
based on  \cite[1.5]{keeler} which yields a stronger conclusion:

\begin{lemma}\label{lemma:keeler}
  Suppose that  $\LL$ is an arithmetically nef line bundle on a projective
variety $X$ with ample line bundle $O(1)$. Then 
there exists a thickening $(\tilde X, \tilde \LL, \tilde O(1))$ such
that the regularities
of $\{\LL_q^m\mid q\in Spec\, A\text{ closed}, m\ge 0 \}$  are bounded above.
\end{lemma}

We will say that the locally free sheaf $\E$ is arithmetically nef if
and only $O_{\PP(E)}(1)$ is. By applying \cite[6.2.12,
6.2.16]{lazarsfeld} to the fibers of a thickening, we obtain
some useful criteria for this condition.

\begin{prop}
  \begin{enumerate}
\item[]
  \item The class of arithmetically nef sheaves is stable under
quotients, extensions and tensor products.
\item If $f:Y\to X$ is a surjective map of projective varieties,
then a locally free sheaf $\E$  on $X$ is arithmetically nef  if and
only if $f^*\E$ is.
\item A locally free sheaf $\E$  is arithmetically nef
if and only if $S^n(\E)$ is for some $n>0$.
  \end{enumerate}
\end{prop}

\begin{cor}
  $\E$ is arithmetically nef if some positive symmetric power is
globally generated.
\end{cor}

Given a locally free sheaf $\E$ and a $\Q$-divisor $D$, it is convenient
to say that a formal twist $\E\langle D\rangle$ \cite[6.2A]{lazarsfeld} is 
arithmetically nef if $S^n(\E)\otimes O(nD)$ is arithmetically
for some $n>0$ such that $nD$ is integral. The choice of $n$ is
immaterial by above.

The key estimate on Frobenius amplitude
is contained in the following theorem.

\begin{thm}\label{thm:key}
Let $X$ be a smooth projective variety defined over a perfect field
$k$ of characteristic $p>0$.
 Let $O_X(1)$ be a sufficiently ample line bundle, and let $\E,\F$ be 
 locally free sheaves on $X$. If $\iota >\lambda(\E(-\dim X))$ and
if $p^N\ge reg(\F)$, then
$$
H^\iota (X,\E^{(p^N)}\otimes \F) = 0
$$
\end{thm}

\begin{proof}
Let $n = \dim X$, and  fix $N>0$ as above. 
Set $f'= F^N$ where $F:X\to X$ is the absolute 
Frobenius (which acts by identity on $X$ as a set, and by $p$th powers
on $O_X$). The endomorphism $f'$ is not $k$-linear, but this can be rectified
by changing the $k$-scheme structure. Let $X'$ be the fibered product 
$X\times_{Spec k}{Spec\, k}$ over the $p^N$th power map of $k$. Then $f'$ factors as
$$X\stackrel{f}{\longrightarrow} X'\stackrel{g}{\longrightarrow} X$$
where $f$ is $k$-linear and $g:X'\to X$ is the natural map.
The morphism $f$ is the relative $p^N$th Frobenius. Note that $f$ is
flat since $X$ is smooth, and that $g$ is an isomorphism of schemes
because $k$ is perfect. Therefore
$g^*,g_*$ induce isomorphisms on cohomology. Set $\E'=g^*\E$ and
 $O_{X'}(1) = g^*O_X(1)$, the latter is easily seen to be  sufficiently ample.

Set
 $$C^i  =
 \begin{cases} \RR_{-i}\boxtimes O_{X'}(i) &\text{if $0\le i>-n-1$}\\
     ker[\RR_n\boxtimes O_{X'}(-n) \to \RR_{n-1}\boxtimes O_{X'}(-n+1) ]
     &\text{if $i=-n-1$}
  \end{cases}
  $$ 
where $\RR_i$ is defined as in section 1. 
These fit into  a resolution $C^\dt$ of the structure sheaf of the diagonal $\Delta$.
In general, for locally free sheaves $\E_i$,
  $(\E_1\boxtimes \E_2)\otimes C^\dt$ is quasiisomorphic to 
$\delta_*(\E_1\otimes \E_2)$,
where $\delta:X'\to X'\times X'$ is the diagonal embedding.
Consequently,
$D^\dt= (O_{X'}(-n)\boxtimes O_{X'}(n))\otimes C^\dt$ gives another
resolution of the diagonal.

 Let $\gamma:X\to X'\times X$ be the
morphism for which $p_1\circ\gamma = f$ and $p_2\circ \gamma = id$, where
$p_i$ are the projections. Let $\Gamma$ be transpose of the graph of $f$,
i.e. the image of $\gamma$ with its
reduced structure. Then $(1\times f)^{-1}\Delta = \Gamma$, and
$(1\times f)^*D^\dt$ gives a resolution for $O_{\Gamma}$ by flatness
of $1\times f$.  Then we have a
quasiisomorphism 
\begin{equation}
  \label{eq:OGamma}
  (\E'\boxtimes \F)\otimes (1\times f)^*D^\dt
\cong 
\gamma_{*}(f^*\E'\otimes \F)= \gamma_{*}(\E^{(p^N)}\otimes \F).
\end{equation}
Thus, we can compute 
\begin{equation}
  \label{eq:EpNXX}
H^i(X,\E^{(p^N)}\otimes \F) =
H^i(X'\times X, \gamma_{*}(f^*\E'\otimes \F))
\end{equation}
using this resolution.

For $a\ge -n$, we have
$$
 (\E'\boxtimes \F)\otimes (1\times f)^*D^a
\cong (\E'\otimes \RR_{-a}(-n))\boxtimes (\F\otimes O_{X}(p^N(n+a)))
$$
 K\"unneth's formula implies
 \begin{equation}
   \label{eq:kunneth}
H^{\iota-a}( \E'\otimes \RR_{-a}(-n)\boxtimes \F(p^N(n+a))) = 
\bigoplus_{r+s=\iota-a} H^r(\E'\otimes \RR_{-a}(-n))\otimes
H^s(\F(p^N(n+a)))   
 \end{equation}
We claim that the above cohomology groups vanish for
$\iota>\lambda(\E(-n))$. We have $a\ge -n$. When $a=-n$,
 cohomology in degree $\iota-a>n$ is automatically zero.
Therefore we can assume $0\ge a> -n$. Now $r+s = \iota-a >
\lambda(\E(-n))$ implies $r>\lambda(\E(-n))$ or  $s>\lambda(\E(-n))\ge
0$. In the first case, the right side of (\ref{eq:kunneth}) vanishes
 by  corollary \ref{cor:lambdaRi}. When  $s>0$, 
the second group on the right vanishes, since $p^N(a_n)\ge reg(\F)$. 
Thus
$$H^{\iota-a}( \E'\otimes \RR_{-a}(-n)\boxtimes \F(p^N(n+a))) = 0$$
as claimed.  

Therefore 
$$
H^\iota(X,\E^{(p^N)}\otimes \F) = 0
$$
for $\iota >\lambda(\E(-n))$ and $p^N\ge reg(\F) $ by
 (\ref{eq:EpNXX}), the above claim, and 
the spectral sequence
$$E_1=H^{\iota-a}( (\E'\boxtimes \F)\otimes (1\times f)^*D^{a}) 
\Rightarrow H^{\iota}(X'\times X, \gamma_{*}(f^*\E'\otimes \F))$$
\end{proof}

\begin{cor}\label{cor:key}
Let $X$ be a smooth projective variety defined over a perfect field
$k$. Let $O_X(1)$ be sufficiently ample, and $\E$ a locally free  sheaf on  $X$.
Then the Frobenius amplitude $\phi(\E)$ satisfies 
 $$\phi(\E)\le \lambda(\E(-\dim X))$$  
\end{cor}

\begin{proof}
  By specialization,  we can assume that  $k$ is a perfect field of
  characteristic $p>0$. The  result is now an immediate consequence of
the theorem.
\end{proof}

Even though the following corollary is special case of the next
theorem, we give a short direct proof.

\begin{cor}
Let $X$ be a smooth projective variety defined over a perfect field
$k$, and let $\LL$ be a line bundle on $X$.
Then   $\phi(\LL)\le \alpha(\LL)$ holds
under one of the following assumptions
\begin{enumerate}
\item $char\, k=p>0$. 
\item $char\, k=0$ and $\LL$ is arithmetically nef.
\end{enumerate}

\end{cor}

\begin{proof}
Let $O_X(1)$ be a sufficiently ample line bundle.
  For some $m_0>0$, we have 
  \begin{equation}
    \label{eq:lambdalealpha}
\lambda(\LL^m(-\dim X)) \le \alpha(\LL^m)\le \alpha(\LL)    
  \end{equation}
whenever $m> m_0$.
The previous corollary yields
$$\phi(\LL^m)\le \lambda(\LL^m(-\dim X))$$
and therefore
$$\phi(\LL^m)\le \alpha(\LL).$$

When  $char\, k=p>0$, we can choose $m$ to be a power $p^N$,
then
$$\phi(\LL^{p^N}) = \phi(\LL)$$
so we are done in this case.

Suppose $char\, k=0$. Choose a thickening $(\tilde X,\ldots)$ of
$(X,\ldots)$ over $Spec\, A$. By semicontinuity, there
is an open set $U\subset Spec\, A$ such that (\ref{eq:lambdalealpha})
holds with $\LL$ replaced by $\LL_q$ on the left, for all closed points
$q\in U$. It would follow that
$$\phi(\LL_q^m)\le \alpha(\LL)$$
Pick $q\in U$, and let $p =p(q)$.
Choose $p^N\ge m$, then by \cite[4.5]{arapura} the $F$-semipositivity of $\LL$ implies
$$\phi(\LL_q)=\phi(\LL_q^{p^N})=\phi(\LL_q^{m}\otimes \LL_q^{p^N-m})\le \phi(\LL_q^{m}) $$
and we are done in this case.
\end{proof}

We now come to our main result, which is a refinement of
\cite[6.1]{arapura}.

\begin{thm}\label{thm:main}
Let $X$ be a smooth projective variety
defined over a   field $k$ of characteristic $0$.
If $\E_1,\ldots \E_m$ are arithmetically nef  vector bundles over $X$, then
  the Frobenius amplitude 
$$\phi(\E_1\otimes \ldots \otimes \E_m) \le
\sum_i(rank(\E_i) + \alpha(\E_i)) -m $$
\end{thm}

We set up the notation for the proof, which  is modeled on the proof of  \cite[6.1]{arapura}. 
Let $r_i= rank(\E_i)$. Choose a sufficiently ample line bundle $O_X(1)$,
and a thickening
$(\tilde X,\tilde \E_1,\ldots \tilde \E_m, \tilde O_X(1))$ of 
$(X, \E_1,\ldots \E_m, O_X(1))$ over
$Spec\, A$. We can assume that the sheaves $\tilde \E_i$ are locally free.
We have the following notation:
 $\BS^\lambda$ is the Schur functor
associated to a partition $\lambda=(\lambda_1,\lambda_2,\ldots)$, 
if $q\in Spec A$, then
 $\E_{h,q}$ is shorthand for 
$ (\tilde \E_h)|_{X_q}$.
 Then \cite[theorem 6.2]{arapura} gives a quasi-isomorphism between
$\E_{h,q}^{(p(q))}$ and a Carter-Lusztig complex
\begin{equation}
  \label{eq:carter}
CL(\E_{h,q})^\dt:= \BS^{(p(q))}(\E_{h,q})\to  
\BS^{(p(q)-1,1)}(\E_{h,q})\to 
\ldots \BS^{(p(q)-r_h+1,1,\ldots 1)}(\E_{h,q}) \to 0
\end{equation}
for any closed point $q\in Spec\, A$ with $p(q)$ large enough.

\begin{lemma}\label{lemma:HiSchur}
 For any integer $m$,
there exists $N_0$ and a nonempty open set $U\subseteq Spec\, A$ such that
$$H^i(X_q, \BS^{(N+1,1,\ldots 1) }(\E_{h,q})\otimes \tilde O_{X_q}(m)) =0 $$
for all $i>\alpha(\E_h)$,  $N\ge N_0$
and $q\in U$.  
\end{lemma}
 
\begin{proof}
The argument is basically a modification of  [loc. cit, 6.5]. 
To avoid excessive notation, we will suppress ``$q$'' and other
decorations whenever it is clear from context what is meant. 
Let $\pi:Flag(\E_h)\to X$ be the flag bundle  of $\E_h$.
Set  $\M= O_X(m)$ and $\LL = O_{\PP(\E_h)}(1)$.
Let  $\kappa=(1\ldots 1)$ be a  partition of length at most $r_h$.
The map $\pi$ factors through a natural map $\pi_1:Flag(\E_h)\to \PP(\E_h)$.
There is a canonical ample  line bundle
$L_\kappa$ on $Flag(\E_h)$ such that $\pi_*L_\kappa  =\BS^\kappa(\E_h) $,
and more generally $\pi_*(\pi_1^*O_{\PP(\E_h)}(N)\otimes L_{\kappa})
=\BS^{(N+1,1,\ldots 1)}(\E_h) $ for $N\ge 0$.
 Kempf's vanishing
theorem  \cite[II, 4.5]{jantzen}, implies that the higher direct
images of $\pi_1^*O_{\PP(\E_h)}(N)\otimes L_{\kappa}$ under $\pi$ and $\pi_1$
vanish. Therefore we have 
\begin{eqnarray*}
 H^i(X_q, \BS^{(N+1,1,\ldots 1) }(\E_{h})\otimes \M) &\cong& 
H^i(Flag(\E_{h,q}), \pi_1^*\LL^N\otimes L_{\kappa}\otimes
\pi^*\M)\\
&\cong& H^i(\PP(\E_{h,q}), \LL^N\otimes\F)\\
&\cong& H^i(\PP(\E_{h,q}), \LL^{N_0}\otimes \LL^{N- N_0}\otimes\F)
\end{eqnarray*}
where $\F=\pi_{1*}(L_k\otimes \M)$ on $\PP(\E_h)$.
The regularities of $(\LL^{N- N_0}\otimes\F)|_{\tilde X_q}$ are bounded above
for $N-N_0\ge 0$ by 
lemma~\ref{lemma:keeler}. Therefore,  by
lemma~\ref{lemma:lambdatensorreg},
 there is $N_0$ such that   the $i$th cohomology of $\LL^{N_0}\otimes
 \LL^{N- N_0}\otimes\F$ vanishes for $N\ge N_0$  and
   $i> \alpha(\E_h)=\alpha(\LL)$.
\end{proof}

\begin{proof}[Proof of theorem~\ref{thm:main}.] 
By the last lemma,
it follows that after shrinking $Spec\, A$, the levels
$$\lambda(\BS^{(p(q)-j, 1,\ldots 1)}(\E_{h,q})(-\dim X)\le \alpha(\E_h)$$ 
for $j=0,\ldots r_h-1$.  Therefore, by corollary \ref{cor:key},
the $F$-amplitudes of the sheaves
$\BS^{(p(q)-j, 1,\ldots 1)}(\E_{h,q})$
occurring in (\ref{eq:carter})
are at most $\alpha(\E_h)$. 
Since $F$-amplitude is subadditive under tensor products [loc. cit, theorem 4.1],
it follows that entries of
\begin{equation}
  \label{eq:tensorcarter}
CL(\E_{1,q})^\dt\otimes \ldots \otimes CL(\E_{m,q})^\dt  
\end{equation}
have $F$-amplitude bounded by $\sum \alpha(\E_h)$. Moreover,  the complex
(\ref{eq:tensorcarter}) is a length $\sum (rank(\E_h)-1)$ resolution
of 
$$(\E_1\otimes \ldots \otimes \E_m)_q^{(p(q))} = \E_{1,q}^{(p(q))}\otimes \ldots 
\otimes \E_{m,q}^{(p(q))}. $$  
Therefore
 [loc. cit, theorem 2.5]
implies the desired inequality 
$$\phi((\E_1\otimes \ldots \E_m)_q) 
=\phi((\E_1\otimes \ldots \E_m)_q^{(p(q))} )
 \le \sum \alpha(\E_h) + \sum (rank(\E_h)-1).$$
\end{proof}

 The following is a 
generalization of \cite[theorem 6.7]{arapura}.
Let us say that  a vector bundle $\E$ is arithmetically nef along a
reduced divisor $D=\sum D_i$
if $\E\langle -\sum r_iD_i \rangle$ is arithmetically
nef for some $0\le r_i< 1$.

\begin{thm}\label{thm:logmain}
Let $X$ be a smooth projective variety 
defined over a   field $k$ of characteristic $0$.
Suppose that $\E$ is a locally free sheaf which is arithmetically
nef along a reduced divisor $D$ with normal crossings.
Then 
$$\phi(\E, D)\le  \alpha(\E,D)+rank(\E)$$
\end{thm}

\begin{proof}
We show that for almost all $q$
there exists a divisor $0\le G\le (p(q)-1) D$ 
such that 
$$\lambda(\BS^{(p(q)-j, 1,\ldots 1)}(\E_{q})(-\dim X)(-G))\le
\alpha(\E, D)$$ 
This is will yield a bound on $F$-amplitudes of the entries of the Carter-Lusztig complex
$CL(\E_{q})^\dt(-G)$
as above.

Let $P = \PP(\E)$, $\LL=O_{P}(1)$ and fix a very ample line
bundle $H$ on $P$. 
By assumption, $\LL^{\ell}(-D')$ are arithmetically nef
for some $0\le D'\le (\ell-1) D$.
We have 
$$\alpha(\E,D) = \alpha(\LL^{M}(-D_h''))$$
for some $M>0$ and $0\le D''\le (M-1)D$.
Fix a constant $c$ that will be determined later.
Then for  $a\gg 0$,
$$H^i(P, \LL^{aM+b}(- aD'')\otimes H^c) = 0$$
for all $i>\alpha(\E,D)$, $b\in\{0,\ldots \ell-1\}$.
The collection of exponents $N_b = aM+b$ is a complete set  of
representatives for $\Z/\ell\Z$. Set $F= aD''$. This
satisfies  $0\le F \le (N_b-1)D $.
We can find a nonempty
open set $U\subset Spec\, A$ such that
$$H^i(P_{q}, \LL_{q}^{N_b}(- F)\otimes H^c) = 0$$
for all $i>\alpha(\E,D)$,  $b$
and closed $q\in U$.

Given $N\ge M = \max N_b$,  choose  $b$
such that $N\equiv N_b\mod \ell$, and set
$$G(N) = F + \frac{N- N_b}{\ell} D'$$
Arguing as in lemma~\ref{lemma:HiSchur}, we obtain that
\begin{equation}
  \label{eq:logmain}
 H^i(X_q, \BS^{(N+1,1,\ldots 1) }(\E)(-G(N))\otimes O(r)) =
H^i(P_q, \LL^{N_b}(-F)\otimes\F)  
\end{equation}
where 
$$\F=\F(N,r) = (\LL^\ell(-D'))^{(N- N_b)/\ell}\otimes \pi_{1*} (L_{(1,\ldots
1)}\otimes \pi^*O(r)) $$
Since the regularities of the sheaves in $\{\F(N,r)_q\mid N\ge M, q\in
U\}$
are bounded above by lemma~\ref{lemma:keeler},
 we can choose the initial constant $c\ll 0$ so that
the cohomologies in  (\ref{eq:logmain}) vanish for
$i>\alpha(\E,D)$ and $r=-2\dim X, \ldots -\dim X$ by 
lemma~\ref{lemma:lambdatensorreg}. This
will force the levels
$$\lambda(\BS^{(p(q)-j, 1,\ldots 1)}(\E_{h,q})(-G(p(q)))(-\dim X))\le
\alpha(\E, D)$$ 
for $p(q)\ge max(M,rank(\E))$ as required.
\end{proof}

The following is a refinement of \cite[theorem 7.1]{arapura}. The hypothesis
about global generation of  $S^{rN}(\E)\otimes \det(\E)^{-N}$ below can
be interpreted as a strong semistability condition, see \cite[p. 247]{arapura1}.

\begin{thm}\label{thm:SrnEdetE}
  Let $\E$ be  a rank $r$ vector bundle on a smooth projective variety $X$
such that $S^{rN}(\E)\otimes \det(\E)^{-N}$ and $\det(\E)^N$ are
globally generated for
some $N> 0$ prime to $char\, k$. Then $\phi(\E)\le \alpha(\det(\E))$.
\end{thm}

\begin{proof}
By specialization, it is enough to  prove the theorem over a 
field with $char\, k=p>0$.
Since $\Sigma_{N}=S^{rN}(\E)\otimes \det(\E)^{-N}$ is globally generated,
$S^{rN}(\E)$ is a quotient of $H^0(\Sigma_N)\otimes \det(\E)^N$. Therefore
$$\alpha(S^{rN}(\E))\le \alpha(H^0(\Sigma_N)\otimes \det(\E)^N)\le \alpha(\det(\E)).$$
The first inequality follows from corollary \ref{cor:quot} and the last
from \cite[cor 1.10]{sommese} (the proof of which is characteristic free).
This implies that $\alpha(\E)\le \alpha(\det(\E))$ by 
corollary \ref{cor:Snample}.

Fix a sufficiently ample line bundle $O_X(1)$, and
choose $q=p^n\equiv 1\, (mod \, N)$. Then as in the proof of 
\cite[theorem 7.1]{arapura}, we can conclude that $\E^{(q)}$ is a direct
summand of $S^q(\E)$. Therefore 
$$H^{i+j}(X,\E^{(q)}\otimes O_X(-\dim X-1-j)) = 0$$ 
for $i> \alpha(\det(\E))$ and $q\gg 0$ subject to the earlier constraint.
In other words, we have shown
 $$\lambda(\E^{(q)}(-\dim X)) \le \alpha(\det\E)$$
By corollary \ref{cor:key}, 
$$\phi(\E)=\phi(\E^{(q)})\le \alpha(\det(\E))$$
\end{proof}

\begin{cor}
 Suppose that $S^{rN}(\E)\otimes \det(\E)^{-N}$ is globally generated
for some $N> 0$ prime to $char\, k$,
and that   some power of $\det(\E)$ is globally generated. 
Then $\E$ is $F$-semipositive.
\end{cor}

\begin{proof}
It is enough to check this in positive characteristic.
Let $\E_n = \E^{(p^n)}\otimes O_X(1)$. Then $\E_n$ satisfies the
hypothesis of the theorem with $\det(\E_n)$ ample for all $n$.
Therefore $\E_n$ is $F$-ample (i. e. has $F$-amplitude $0$).
This implies that the  regularities of $\E^{(p^n)}$ are bounded away from
$\infty$. 
\end{proof}

\section{Vanishing theorems}\label{section:vanish}
 
In this section, we fix a smooth $n$ dimensional projective
variety $X$ defined  over a field $k$ of characteristic $0$.
Applying \cite[cor 8.6]{arapura} to theorems \ref{thm:main},
\ref{thm:SrnEdetE}, and \ref{thm:logmain} yields:

\begin{thm}\label{thm:van1}
  Suppose that
\begin{enumerate}
\item  $\E_1,\ldots \E_m, \F$ are locally free sheaves on $X$ or ranks
$r_1,\ldots r_m,s$ respectively,
\item $\E_1$ is arithmetically nef along a reduced divisor with normal
   crossings $D$,
\item $S^{sN}(\F)\otimes \det(\F)^{-N}$ and $\det(\F)^N$ are
globally generated for some $N>0$. 
\end{enumerate}
Then 
$$H^q(X,\Omega_X^p(\log D)(-D)\otimes \E_1\otimes\ldots\otimes
\E_m\otimes \F) = 0$$
for 
$$p+q> n+\alpha(\E_1,D) + \alpha(\E_2)+\ldots \alpha(\E_m) +\sum r_i -m.$$
\end{thm}

\begin{cor}\label{cor:van1}
Suppose that  $D= D_1+\ldots D_d$ is a divisor with normal crossings
and $\E$ a vector bundle such that $\E$ is  arithmetically nef
along $D$, and each $\E|_{D_i}$ is arithmetically nef along
$D_i'= D_{i+1}+\ldots D_d$. Then
 $$H^q(X,\Omega_X^p\otimes \E) = 0$$
for 
$$p+q\ge n+ rank(\E)+max(\alpha(\E,D),\alpha(\E|_{D_1},D_1')-1,
\alpha(\E|_{D_2}, D_2')-1,\ldots);$$ 
in particular for
$$p+q\ge n+rank(\E) + \alpha(\E)$$
\end{cor}

\begin{proof}
For the last part, we note the inequality
$$max(\alpha(\E,D),\alpha(\E|_{D_1},D_1')-1,\ldots)\le \alpha(\E).$$
Tensoring $\E$ with   the sequences
$$0\to \Omega^p_X(\log D_{i-1}')(-D_{i-1}')\to \Omega^p_{X}(\log D_i')(-D_i')\to 
 \Omega^p_{D_i}(\log D_i')(-D_i')\to 0$$
and applying the theorem implies the vanishing statement.
\end{proof}

The special case of the above result, when $\E$ is $d$-ample (and $D=0$), is
due to Sommese and Shiffman \cite[cor. 5.20]{ss}.

\begin{cor}\label{cor:vanalphagen}
Suppose that  $\E$ is a vector bundle such that its
 pullback under a birational map $f:Y\to X$, with $Y$ smooth,  is  arithmetically nef
 along a normal crossing divisor. Then
we have
 $$H^i(X,\omega_X\otimes \E) = 0$$
if $i\ge \alpha_{gen}(\E) + rank(\E)$, where $\omega_X= \Omega_X^n$.  
\end{cor}

Using the method of Manivel, we can get a vanishing theorem for tensor powers.
We start with a preliminary lemma.

\begin{lemma}\label{lemma:proofofvan2}
  Let $\E$ and $\F$ be a vector bundles on $X$ with $\F$ 
$F$-semipositive. Let $P$ be the $m$-fold fiber product 
$\PP(\E)\times_X\ldots \times_X\PP(\E)$ with projections
$p_i:P\to \PP(\E)$ and $\pi:P\to X$. Then for any sequence of positive integers
$r,j_1,\ldots j_m$,
$$H^q(P,\Omega^p_P\otimes p_1^*O(j_1)\otimes \ldots
p_m^*O(j_m)\otimes \pi^*\det(\E)^{\otimes r}\otimes \pi^*\F)=0 
$$
for $p+q> \dim P+\alpha_\otimes(\E)=m(rank(\E)-1)+\dim X +\alpha_\otimes(\E)$.
\end{lemma}

\begin{proof}
  Let 
$$\G = S^{j_1}(\E)\otimes \ldots S^{j_m}(\E)\otimes (\det \E)^{\otimes r}.$$
Then since $\G$ is a summand of a tensor power of $\E$, we have
$\alpha_{\otimes}(\G) \le \alpha_{\otimes}(\E)$.
Furthermore,
$\alpha(O_{\PP(\G)}(1))=\alpha(\G)\le \alpha_{\otimes}(\G)$ by 
proposition~\ref{prop:alphaVsalphaOtimes}.

We have an embedding $P\to \PP(\G)$ such that
$\alpha_{\PP(\G)}(1)$ restricts to  $p_1^*O(j_1)\otimes \ldots
p_m^*O(j_m)\otimes \pi^*\det(\E)^{\otimes r}$.
Therefore by corollary~\ref{cor:restalpha} 
$$\alpha( p_1^*O(j_1)\otimes \ldots \otimes \pi^*\det(\E)^{\otimes r})
\le \alpha_{\otimes}(\E)$$
So by  \cite[prop 3.10, theorem 4.5]{arapura} and theorem~\ref{thm:main}
 $$\phi(p_1^*O(j_1)\otimes \ldots \otimes \pi^*\det(\E)^{\otimes r}\otimes\pi^*\F)
\le \alpha_{\otimes}(\E)$$
Therefore  \cite[cor 8.6]{arapura} implies that 
$$H^q(P,\Omega^p_P\otimes p_1^*O(j_1)\otimes \ldots
p_m^*O(j_m)\otimes \pi^*\det(\E)^{\otimes r}\otimes \pi^*\F)=0 
$$
vanishes for 
$$p+q > \dim P + \alpha_{\otimes}(\E)\ge \alpha( p_1^*O(j_1)\otimes \ldots
p_m^*O(j_m)\otimes \pi^*\det(\E)^{\otimes r})$$

\end{proof}

\begin{thm}\label{thm:van2}
  Let $\E$ and $\F$ be a vector bundles on $X$.
Suppose that $\F$ is $F$-semipositive. 
Then for any sequence of integers $k_1,\ldots k_\ell, j_1,\ldots j_m$
$$H^q(X,\Omega_X^p\otimes S^{k_1}(\E)\otimes \ldots S^{k_\ell}(\E)\otimes
\wedge^{j_1}\E\otimes\ldots \wedge^{j_m}\E\otimes (\det
\E)^{\ell+n-p}\otimes \F) = 0$$
whenever
$$p+q> n +\sum(rank(\E) - j_i) + \alpha_\otimes(\E).$$
\end{thm}

\begin{proof}

The proof of theorem $A$  given in \cite[section 2.2]{manivel}
 can be carried out  almost word for
word in the present context,  with $\LL$ replaced by $\F$
and the appeal to  Kodaira-Akizuki-Nakano by the use of  
lemma \ref{lemma:proofofvan2}.

\end{proof}

\begin{cor}
  If $\E$ is globally generated and $ \F$ $F$-semipositive,
then 
$$H^q(X,\Omega_X^p\otimes S^{k_1}(\E)\otimes \ldots S^{k_\ell}(\E)\otimes
\wedge^{j_1}\E\otimes\ldots \wedge^{j_m}\E\otimes (\det
\E)^{\ell+n-p}\otimes \F) = 0$$
whenever 
$$p+q> n +\sum(rank(\E) - j_i)+ \alpha(\E) $$
\end{cor}

\begin{proof}
  This follows from proposition \ref{prop:alphaVsalphaOtimes}.
\end{proof}

We can extend the above results to certain geometric variations of 
Hodge structure.
Recall that a morphism $f:X\to Y$ of smooth projective
varieties is called semistable if it is given
local analytically (or \'etale locally)  by
$$y_1 = x_1x_2\ldots x_{d_1},\> y_2= x_{d_{1}+1}x_{d_1+2}\ldots x_{d_{2}},\>\ldots$$
It follows that the discriminant $D\subset Y$ and its preimage $E\subset X$
are reduced divisors with normal crossings. Let 
$\Omega_{X/Y}^\dt(\log E/D)$ denote the relative logarithmic de Rham complex.
Then
$$H = \R^h f_*\Omega_{X/Y}^\dt(\log E/D)$$
forms part of a variation of Hodge structure. We have a Gauss-Manin connection
$$\nabla :H\to \Omega_Y^1(\log D)\otimes H,$$
and a Hodge filtration
$$F^aH = image[\R^h f_*\Omega_{X/Y}^{\ge a}(\log E/D)\to H],$$
satisfying Griffiths' transverality 
$$\nabla(F^aH)\subset \Omega_Y^1(\log D)\otimes F^{a-1}H.$$
The connection $\nabla$ is integrable, so it fits into a complex
$$\Omega^\dt( H) =H\to \Omega_Y^1(\log D)\otimes H\to  
\Omega_Y^2(\log D)\otimes H\ldots $$
which, by Griffiths' transversality, is compatible with the filtration
when suitably shifted.
The associated graded complex $Gr^a(\Omega^\dt(H))$ can be identified with
the complex 
 $$
 \Omega_Y^\dt(\log D)\otimes R^{h+\dt}f_*\Omega_{X/Y}^{a-\dt}(\log E/D),
$$
where the differentials are given by cup product with the Kodaira-Spencer
class of $f$.
Using the same method as in the proof of \cite[cor. 8.6]{arapura}
but replacing \cite{deligne-ill} by \cite[theorem 4.7]{illusie},
we obtain:

\begin{thm}
With the above assumptions and notation,
 for any coherent sheaf $\E$ on $Y$,
$$\HH^{b+a}(Y,Gr^a(\Omega^\dt(H))\otimes \E) = 0 $$
for $a+b \ge \dim Y + \phi(\E)$.
\end{thm}

\begin{cor}
If $\E$ is locally free, then
  $$\HH^{b+a}(Y,Gr^a(\Omega^\dt(H))\otimes \E) = 0 $$
for $a+b \ge \dim Y + rank(\E) +\alpha(\E)$.
\end{cor}

Taking $a$ to be maximal yields the  following
 Koll\'ar type vanishing theorem:

\begin{cor}\label{cor:kollar}
  If $\E$ is locally free, then
$$H^i(Y, R^hf_*\omega_X\otimes \E) = 0$$
for $i\ge rank(\E) + \alpha(\E)$.
\end{cor}

We conjecture that the last corollary
 should hold without the semistability hypothesis
on $f$. As an immediate consequence of the
semistable reduction theorem \cite{kkms}, we see
that this is so when $Y$ is a curve:

\begin{cor}
  Corollary  \ref{cor:kollar} holds for an arbitrary
map when $\dim Y =1$.
\end{cor}

 We also conjecture that corollary \ref{cor:kollar} holds when
$\alpha(\E)$ is replaced by $\alpha_{gen}$. The special
case when $f=id$
is proved in corollary \ref{cor:vanalphagen}.

\section{Varieties with endomorphisms}

 Fix an integer $q$.
By a {\em $q$-endomorphism} of a smooth variety $X$, 
we mean a  morphism $\phi:X\to X$ such that
$\phi^*:NS(X)\otimes \Q\to NS(X)\otimes \Q$ acts by multiplication
by $q$, where $NS(X)$ is the Neron-Severi group.

\begin{ex}\label{ex:qmorphPn}
  The morphism $\PP^n\to \PP^n$  given by 
$$[x_0,\ldots x_n]\mapsto [x_0^q,\ldots x_n^q]$$
is a $q$-endomorphism.
\end{ex}

\begin{ex}
    Let $X$ be an Abelian variety, and let $\mu_{q,X}:X\to X$ be the 
 multiplication by $q$. This is a $q^2$-endomorphism \cite[p. 75]{mumford}.
\end{ex}

\begin{ex}
  The product of two varieties with $q$-endomorphisms
  has a $q$-endomorphism.
\end{ex}

Before giving the  final example,  we remind the reader that there is
a recipe for  building  a toric variety
from a fan $\Delta$   in a lattice $N$ \cite{danilov,fulton, oda}.
This is denoted by $T_Nemb(\Delta)$ or $X(\Delta)$.
This has an action of the torus $T_N = Hom(\check{N},\C^*)$.
A $T_N$-equivariant  divisor $D$ is determined 
by a piecewise linear real valued function  on the support of $\Delta$
taking integer values on $N$ \cite[sect 2.1]{oda}.
 The collection of these forms an Abelian  group $SF(\Delta)$. 
 Moreover, when $X(\Delta)$ is complete, we have
 a surjective homomorphism $SF(\Delta)\to Pic(X(\Delta))$ 
 \cite[cor. 2.5]{oda}.
 
\begin{ex}\label{ex:toric}
    Assume that $X=X(\Delta)$ is smooth and projective.
For any $q$, multiplication by $q$ on $N$ induces a morphism
of fans $\Delta \to \Delta$ \cite[sect. 1.5]{oda}. This induces an
endomorphism of the associated toric variety $\mu_{q,X}:X \to X$.
It is clear that $\mu_{q,X}$ acts by multiplication by $q$ on
$SF(\Delta)$, and this implies that $\mu_{q,X}$ is 
$q$-endomorphism.
For the special case of $\PP^n$ with its standard toric structure,
$\mu_{q,X}$ coincides with the morphism $f$ given in example \ref{ex:qmorphPn}.
\end{ex}

For line bundles theorem \ref{thm:main} extends to $q$-morphisms.

\begin{lemma}
Suppose that $f:X\to X$ is a $q$-morphism of a projective
variety with $q>1$. 
For any line bundle $\LL$ on $X$, $\alpha_f(\LL)\le \alpha(\LL)$ 
\end{lemma}

\begin{proof}
The set of numerically trivial line bundles $Pic^\tau(X)$ forms
a bounded family.
Therefore given a vector bundle $\E$, we can choose $N_0>0$ so that
  $H^i(X,\LL^{\otimes q^N}\otimes \E\otimes \M) = 0$
for $N\ge N_0$, $i> \alpha(\LL)$ and $\M\in Pic^\tau(X)$.  
By assumption, $(f^N)^*\LL = \LL^{\otimes q^N}\otimes \M$ for some $\M\in
Pic^\tau(X)$. Thus $H^i(X,(f^N)^*\LL\otimes \E) = 0$
for $i> \alpha(\LL)$ .
\end{proof}

\begin{thm}\label{thm:sec}
Suppose that $X$ is a smooth projective $k$-variety 
with a $q$-endomorphism $f:X\to X$, such that  $q>1$ is prime to
the characteristic of $k$. Then for any  vector bundle $\E$ on $X$,
$$H^i(X,\Omega_X^j\otimes \E) = 0$$
whenever $i>\alpha_f(\E)$.
\end{thm}

\begin{cor}
    Then for any  line bundle
$$H^i(X,\Omega_X^j\otimes \LL) = 0$$
for $i> \alpha(\LL)$. In particular, this holds for $i>0$
if $\LL$ is ample. 
\end{cor}

For toric varieties, the last part of the  corollary, where is $\LL$ ample,
is due to Danilov \cite[7.5.2]{danilov} (and Bott for $\PP^n$).
For Abelian varieties, where
$\Omega_X^1$ is trivial, the last part is also  well known
\cite[III-16]{mumford}.

Theorem \ref{thm:sec} is an immediate consequence of the
next two lemmas. The first lemma is just a formulation of a well known trick
used in Frobenius splitting arguments \cite{mehta}.

    We call a coherent sheaf $\F$ {\em $f$-split}, if there exists
 a split injection $\F\to f_*\F$.

\begin{lemma}\label{lemma:split}
 If $\F$ is    $f$-split, then
 $H^i(X,\F\otimes \E) = 0$ for $i>\alpha_f(\E)$.
\end{lemma}   

\begin{proof}
Using the projection formula,  we have an injection
 $$H^i(X,\F\otimes \E)\to  H^i(X,(f^N)_*\F\otimes \E) =
H^i(X, \F\otimes (f^N)^*\E)$$
\end{proof}

\begin{lemma}
Suppose that $X$ is a smooth projective $k$-variety 
with  $q$-endomorphism $f:X\to X$, with $q>1$ prime to
the characteristic of $k$. Then $\Omega_X^j$ is $f$-split.  
\end{lemma}

\begin{proof}
  A splitting of  the natural map
$\iota:\Omega_X^j\to f_*\Omega_X^j$ is
given by the $(1/q^{\dim X})Tr$, where $Tr$ is
Grothendieck's trace \cite[4]{deligneL}. The
identity $Tr\circ\iota = q^{\dim X}$ can checked on open
set $U\subset X$ over which $f$ is etale.
\end{proof}

This technique can be adapted to prove other results.
As an example, we give a relative vanishing theorem.

\begin{thm}
Suppose that $q$ is an integer relatively prime to
the characteristic of $k$.  Let
$$
\xymatrix{
 X\ar[r]^{f}\ar[d]^{\pi} & X\ar[d]^{\pi} \\
 Y\ar[r]^{\phi} & Y
}
$$
be a commutative diagram of projective varieties with
$X$ smooth and $f$ and $\phi$ $q$-morphisms.
 Then for any  vector bundle $\E$ on $Y$
$$H^i(Y,R^\ell\pi_*\Omega_X^j\otimes \E) = 0$$
whenever $i>\alpha_f(\E)$. In particular, this
holds for $i>0$ if $\E$ is an ample line bundle.
\end{thm}

\begin{proof}
The $f$-splitting of $\Omega_X^j$ induces a $\phi$-splitting
on $R^\ell \pi_*\Omega_X^j$.  
\end{proof}

\begin{cor}
The conclusion of the theorem holds when $\pi:X\to Y$ is
a toric morphism of toric varieties.
\end{cor}

\end{document}